\documentstyle[12pt]{article}
\def\beq{\begin{equation}}
\def\eeq{\end{equation}}
\def\bea{\begin{eqnarray}}
\def\eea{\end{eqnarray}}
\def\nn{\nonumber}
\def\bra#1{\left\langle #1\right|}
\def\ket#1{\left| #1\right\rangle}
\def\ket#1{\left| #1\right\rangle}
\def\braket#1#2{\VEV{#1 | #2}}
\def\VEV#1{\left\langle #1\right\rangle}

\def\slt{{\cal U}_h(sl(2))}
\newtheorem{definition}{Definition}

\newtheorem{thm}{Theorem}
\begin{document}

\thispagestyle{empty}
%\begin{flushright}
%OWUAM-024 \\
%November 13, 1997 \\
%revised February 23, 1998
%revised March 27, 1998
%\end{flushright}

\vspace*{2cm}

\begin{center}
{\Large Tensor Operators for $ \slt $}

\vspace{2cm} 

N. Aizawa

\medskip

{\em Department of Applied Mathematics, 

Osaka Women's University, Sakai, Osaka 590-0035, JAPAN}

\end{center}

\vfill 

\begin{abstract}
Tensor operators for the Jordanian quantum algebra $ \slt $ are considered. 
Some explicit examples of them, which are obtained in the boson or fermion 
realization, are given and their properties are studied. 
It is also shown that the Wigner-Eckart's theorem can be extended to $ \slt $. 
\end{abstract}

\newpage
%
%%%%%%%%%%%%%%%%%%%%%%%%%%%%%%%%%%%%%%%%%%%%%%%%%%%%%%%%%
%
%                      Introduction
%
%%%%%%%%%%%%%%%%%%%%%%%%%%%%%%%%%%%%%%%%%%%%%%%%%%%%%%%%%
%
\section{Introduction}

  Recent studies on quantum matrices in two dimensions show that 
the Lie group $ GL(2) $ admits two kinds of quantum deformation 
\cite{dmmz,eow,zak}. One of them is denoted by $ GL_{p,q}(2) $ and 
has been studied extensively since the beginning of quantum group 
theory. The other is denoted by $ GL_{g,h} (2) $, it is  sometimes 
called the Jordanian quantum group. The $ SL_h (2) $ is the 
special case of $ GL_{g,h}(2) $ obtained by setting $ g = h $ and 
the quantum determinant unity. The dual of $ SL_h(2) $ is  a 
deformation of the universal enveloping algebra of $ sl(2) $ \cite{ohn} 
and it is called the Jordanian quantum algebra $ \slt $. The explicit form 
of universal $R$-matrix for $ \slt $ is known \cite{bh,sak}. 
It is also known that the $ \slt $ can be obtained from the Drinfeld-Jimbo's 
$ {\cal U}_q (sl(2)) $ by a contraction \cite{aks}. The Hopf algebra dual to 
the $ GL_{g,h} (2) $ was found very recently \cite{adm}.

  The representation theory of $ \slt $ seems to have attracted some interest, 
since it has been revealed that the representation theories of $ \slt $ and 
$ sl(2) $ have some similarities.  
The finite dimensional irreducible representations (irreps) was first considered 
in \cite{dov}, then a simple way to construct irreps with a nonlinear relation 
between the generators of $ \slt $ and $ sl(2) $ was proposed \cite{acc1}. 
They show that the finite dimensional irreps of $ \slt $ can be classified in the same 
way as those of $ sl(2) $ (see also \cite{acc2}). The infinite dimensional representations 
are considered in \cite{bhn} with boson realizations. The first attempt to decompose 
a tensor product of two finite dimensional irreps was made in \cite{naru}, then 
the problem was completely solved in \cite{vdj1,vdj2}. This gives the second similarity 
between the representation theories of $ \slt $ and $ sl(2) $ , that is, the decomposition 
rule is exactly the same as $ sl(2) $. Furthermore a explicit formula for 
$ \slt $ Clebsch-Gordan coefficients (CGC) is given in \cite{vdj1}. 

  The nonlinear relation introduced  in \cite{acc1} gives an interesting observation 
for the coproduct of $ \slt $. We can regard $ \slt $ as the angular momentum 
algebra with a nonstandard coupling rule. This might suggest that $ \slt $ has 
lots of application to various fields in physics. 

   In this paper, we further develop the representation theory of $ \slt $, especially 
tensor operators will be studied. We review the known results on the representation 
of $ \slt $ in the next two sections, in order to fix our notations and to list formulae used in the subsequent 
sections. Tensor operators for $ \slt $ are introduced in \S 4 according to 
\cite{rs}. Some explicit examples of $ \slt $ tensor operators are given and 
their properties are considered. In \S 5, we consider a extension of the Wigner-Eckart's 
theorem to $ \slt $.

%
%%%%%%%%%%%%%%%%%%%%%%%%%%%%%%%%%%%%%%%%%%%%%%%%%%%%%%%%%
%
%               Uh(sl(2)) and its Representations
%
%%%%%%%%%%%%%%%%%%%%%%%%%%%%%%%%%%%%%%%%%%%%%%%%%%%%%%%%%
%
\section{$ \slt $ and Its Representations}

  The Jordanian quantum algebra $ \slt $ is an associative algebra 
with unity and generated by $ X, Y $ and $ H $ subjected to the 
relations \cite{ohn}
\bea
 & & [X,\; Y] = H, \qquad
       [H, \; X] = 2{\sinh hX \over h},
       \nn \\
 & &  [H, \; Y] = -Y (\cosh hX) - (\cosh hX) Y,  \label{define}
\eea
where $h$ is the deformation parameter. The coproduct $ \Delta $, 
the counit $ \epsilon $ and the antipode $ S $ are given by 
\bea
& & \Delta(X) = X \otimes 1 + 1 \otimes X, \nn \\
& & \Delta(Y) = Y \otimes e^{hX} + e^{-hX} \otimes Y, \label{coproduct} \\
& & \Delta(H) = H \otimes e^{hX} + e^{-hX} \otimes H, \nn \\
& & \epsilon(X) = \epsilon(Y) = \epsilon(H) = 0, \label{counit} \\
& &  S(X) = -X, \quad  S(Y) = -e^{hX} Y e^{-hX}, \quad S(H) = -e^{hX}He^{-hX}, 
       \label{antipode}
\eea
so that $ \slt $ is a Hopf algebra. The Casimir element belonging to the 
center of $ \slt $ is \cite{bh}
\beq
  C = {1 \over 2 h} \{ Y(\sinh hX) + (\sinh hX)Y \} + {1 \over 4} H^2 + 
        {1 \over 4} ( \sinh hX )^2. 
  \label{casimir}
\eeq
The Jordanian quantum algebra $ \slt $ is a non-standard deformation of the 
universal enveloping algebra of $ sl(2) $, since all expressions given in (\ref{define})
-(\ref{casimir}) reduce to the corresponding ones for $ sl(2) $ in the limit of 
$ h \rightarrow 0. $ 

  Note that we can eliminate the deformation parameter $h$ from all expressions by 
making the replacement $ hX \rightarrow X $ and $ h^{-1}Y \rightarrow Y $. 
Thus $ \slt $ is isomorphic to $ {\cal U}_{h=1} (sl(2)). $ We however remain the 
parameter $ h $ throughout this paper in order to consider the limit of 
$ h \rightarrow 0 $. 

  The finite dimensional irreps of $ \slt $ can be obtained by using the 
nonlinear relation between generators of $ \slt $ and $ sl(2) $  \cite{acc1}. 
With the definition 
\beq
Z_+ = {2 \over h} \tanh {hX \over 2}, \qquad 
Z_- =  (\cosh {hX \over 2} ) Y ( \cosh{hX \over 2}) ,
\label{sltwo}
\eeq
it follows that $ Z_{\pm} $ and $ H $ satisfy the $ sl(2) $ commutation relations
\beq
   [H, \; Z_{\pm}] = \pm 2Z_{\pm}, \qquad
   [Z_+,\; Z_-] = H,                     \label{sltwo2}
\eeq
and the Casimir element (\ref{casimir}) is rewritten as
\beq
  C = Z_+Z_- + {\displaystyle {H \over 2} \left( {H \over 2} -1 \right)}_.
  \label{casimir2}
\eeq
We can  take undeformed representation bases for $ Z_{\pm} $ and $ H $
\bea
 & & Z_{\pm} \ket{jm} = \sqrt{(j\mp m)(j\pm m+1)} \ket{j m\pm 1}, \nn \\
 & & H \ket{jm} = 2m \ket{jm}, \label{rep1} \\
 & & C \ket{jm} = j(j+1) \ket{jm}, \nn
\eea
where $ j= 0, 1/2, 1, 3/2, \cdots $ and $ m = -j, -j+1, \cdots , j. $ 
The vectors $ \{ \ket{jm} \} $ are nothing but the representation bases for 
$ sl(2) $, their complete orthonormality and the representation matrices 
for bra vectors follow immediately. The representation matrices for $ X $ and 
$ Y $ can be obtained by solving (\ref{rep1}) with respect to $ Z_{\pm}. $ 
The closed form of their expressions is given in \cite{vdj1} and this shows 
that finite dimensional highest weight irreps for $ \slt $ are classified in the 
same way as $ sl(2) $.

%%%%%%%%%%%%%%%%%%%%%%%%%%%%%%%%%%%%%%%%%%%%%%%%%%%%%%%%%%%%
%
%            Clebsch-Gordan Coefficients for Uh(sl(2))
%
%%%%%%%%%%%%%%%%%%%%%%%%%%%%%%%%%%%%%%%%%%%%%%%%%%%%%%%%%%%%
%
\setcounter{equation}{0}
\section{Clebsch-Gordan Coefficients for $ \slt $}

  In this section, we review some known results on the tensor products of 
two irreps given in the previous section. Although $ Z_{\pm} $ 
and $ H $ are the elements of $ sl(2) $, their coproducts are given 
in terms of $ \Delta(X), \Delta(Y) $ and $ \Delta(H) $ 
(see \cite{naru,vdj1} for explicit formulae of $ \Delta(Z_{\pm} $ )) so that 
the irreducible decomposition of tensor product representations 
is a nontrivial problem. 
This problem is solved in \cite{naru,vdj1} and \cite{vdj2}. 

\begin{thm}
  The tensor product of two irreps of $ \slt $ with highest weight 
$ j_1 $ and $ j_2 $ is completely reducible and the decomposition 
into irreps is given by
\beq
  j_1 \otimes j_2 = j_1 + j_2 \oplus j_1 + j_2 - 1 \oplus \cdots \oplus 
  | j_1 - j_2 |, \label{rule}
\eeq
where each irrep is multiplicity free. Namely, the decomposition rules 
for $ \slt $ and $ sl(2) $ are same. 
\end{thm}

  The CGC for $ \slt $ can be obtained by introducing new vectors 
defined by
\beq
  \ket{(j_1m_1) (j_2m_2)} = \sum_{k_i=m_i}^{j_i} \alpha_{k_1, k_2}^{m_1, m_2} 
      \ket{j_1 k_1} \otimes \ket{j_2 k_2}, \label{inter}
\eeq
where the coefficients $  \alpha_{k_1, k_2}^{m_1, m_2} $ are given by
\beq
  \alpha_{k_1, k_2}^{m_1, m_2} = 
      (-1)^{k_2-m_2} {\displaystyle \left({h \over 2}\right)^{k_1+k_2-m_1-m_2} }
      D_{k_1, k_2}^{m_1, m_2} 
      ( b_{k_1, k_2}^{m_1, m_2} - b_{k_1-1, k_2-1}^{m_1, m_2}),
      \label{alpha}
\eeq
with 
\bea
& & D_{k1, k_2}^{m_1, m_2} = 
     {\displaystyle \left\{ 
      { (j_1-m_1)! (j_1+k_1)! (j_2-m_2)! (j_2+k_2)! \over 
        (j_1+m_1)! (j_1-k_1)! (j_2+m_2)! (j_2-k_2)!}
      \right\}^{1/2}_,} \nn \\
& & b_{k1, k_2}^{m_1, m_2} = {\displaystyle 
      \left(
      \begin{array}{c}
      m_1+k_1 \\ k_2-m_2
      \end{array}\right) \; 
      \left(
      \begin{array}{c}
      m_2+k_2 \\ k_1-m_1 
      \end{array}\right)_. } \nn
\eea
We use the following definition of the binomial coefficients, since 
the suffices $ m_i $ in $  b_{k1, k_2}^{m_1, m_2} $ takes  negative values
\[
   \left(
    \begin{array}{c}
    n \\ m
    \end{array}
  \right)
  = 
  \left\{
  \begin{array}{cc}
  \displaystyle{ n(n-1)(n-2) \cdots (n-m+1) \over m!} & 
  {\rm for} \quad m \geq 0 \\
  0 & {\rm for} \quad m < 0
  \end{array}
  \right.
\]
Note that the coefficients $ \alpha_{k_1, k_2}^{m_1, m_2} $ depend on 
$ j_1 $ and $ j_2 $, although the dependence is not shown explicitly. 
Note also that, in the limit of $ h \rightarrow 0 $, all coefficients 
$ \alpha_{k_1, k_2}^{m_1, m_2} $ vanish except for 
$ \alpha_{m_1, m_2}^{m_1, m_2} =1 $ so that 
$ \ket{(j_1m_1) (j_2m_2)} \rightarrow \ket{j_1 m_1} \otimes \ket{j_2 m_2} $. 
We refer to the vectors (\ref{inter}) as "intermediate vectors " in 
this paper, since they appear the intermediate step to the CGC.

  The important property of the intermediate vectors, which plays 
a crucial role in the following discussion, is the action of 
$ \Delta(Z_{\pm}) $ and $ \Delta(H) $ on the intermediate vectors 
given by
\bea
 \Delta(H) \ket{(j_1m_1) (j_2m_2)} &=& 2(m_1+m_2) \ket{(j_1m_1) (j_2m_2)}, 
       \nn \\
 \Delta(Z_{\pm}) \ket{(j_1m_1) (j_2m_2)} 
 &=& \sqrt{(j_1\mp m_1) (j_1 \pm m_1 + 1)} \ket{(j_1\; m_1\pm 1) (j_2m_2)},  
 \label{action} \\
 &+& \sqrt{(j_1 \mp m_2) (j_2 \pm m_2 +1)} \ket{(j_1m_1) (j_2\; m_2\pm 1)}. \nn
\eea
This tells us that the action of $ \Delta(Z_{\pm}) $ and $ \Delta(H) $ on a  
intermediate vector is the same as the action of the undeformed coproducts 
of $ sl(2) $ elements on a vector $ \ket{j_1 m_1} \otimes \ket{j_2 m_2} $. 
Therefore, the bases of irreps for $ \Delta(Z_{\pm}) $ and $ \Delta(H) $ 
are obtained by linear combinations of the intermediate vectors with 
the CGC for $ sl(2) $.
\bea
\ket{jm} &=& \sum_{m_1+m_2=m}\; C_{m_1, m_2, m}^{j_1, j_2,  j} 
  \ket{(j_1m_1) (j_2m_2)}, \nn \\
  &=& \sum_{k_i=m_i}^{j_i}\;  \sum_{m_1+m_2=m}\; C_{m_1, m_2, m}^{j_1, j_2,  j} \; 
         \alpha_{k_1, k_2}^{m_1, m_2} \ket{j_1 k_1} \otimes \ket{j_2 k_2}, \label{bases}
\eea
where $ C_{m_1, m_2, m}^{j_1, j_2,  j} $ is a $ sl(2) $ CGC.

 The orthogonality of the coefficients $ \alpha_{k_1, k_2}^{m_1, m_2} $ is obtained 
in \cite{vdj2} 
\beq
   \sum_{k_1, k_2} \; \alpha_{k_1, k_2}^{m_1, m_2} \alpha_{-k_1, -k_2}^{-n_1, -n_2} 
   = \delta_{m_1, n_1} \delta_{m_2, n_2}.
   \label{orth}
\eeq

  Before closing this section, we add a new result. The intermediate vectors for 
the dual space (space spanned by bra vectors) are given by
\beq
  \bra{(j_1 m_1) (j_2 m_2)} = \sum_{k_i = -j_i}^{m_i} \; 
  \alpha_{-k_1, -k_2}^{-m_1, -m_2} \; \bra{j_1 k_1} \otimes \bra{j_2 k_2}. 
  \label{bra}
\eeq
The action of $ \Delta(Z_{\pm}) $ and $ \Delta(H) $ on (\ref{bra}) 
is the same as the action of the undeformed coproducts 
of $ sl(2) $ elements on a vector $ \bra{j_1 m_1} \otimes \bra{j_2 m_2} $. 
This can be proved by the same way as in \cite{vdj1}. The orthogonality 
of the coefficients $ \alpha_{k_1, k_2}^{m_1, m_2} $ results the orthonormality 
of the intermediate vectors
\beq
   \VEV{(j_1 n_1) (j_2 n_2) | (j_1 m_1) (j_2 m_2)} = \delta_{n_1, m_1} \delta_{n_2, m_2}.
   \label{orth2}
\eeq
Note that the representations of $ \Delta(H) $ and $ \Delta(Z_{\pm}) $ on 
the intermediate vectors (for both bra and ket vectors) are unitary. Therefore 
Eq. (\ref{orth2}) is nothing but the well-known fact that the eigenvectors of 
a hermitian operator with different eigenvalues are orthogonal each other.

%%%%%%%%%%%%%%%%%%%%%%%%%%%%%%%%%%%%%%%%%%%%%%%%%%%%%%%%%%%
%
%              Tensor operators for Uh(sl(2))
%
%%%%%%%%%%%%%%%%%%%%%%%%%%%%%%%%%%%%%%%%%%%%%%%%%%%%%%%%%%%%
%
\setcounter{equation}{0}
\section{Tensor Operators for $ \slt $}
\subsection{Definition of Tensor Operators}

   Rittenberg and Scheunert gave a general definition of tensor operators 
for a Hopf algebra \cite{rs}. To define tensor operators, we first define 
the adjoint action of a Hopf algebra. 
\begin{definition}
Let $ {\cal H} $ be a Hopf algebra, 
let $ W, W' $ be its representation space, and let $ t $ be a operator 
which carries $ W $ into $ W' $. Then the adjoint action of $ c \in {\cal H} $ 
on $ t $ is defined by
\beq
   ad \; c(t) = \sum_i\; c_i \, t \, S(c'_i), \label{ad}
\eeq
where the coproduct for $ c $ is written as 
$ \Delta(c) = \sum_i\; c_i \otimes c'_i $. 
\end{definition}
The adjoint action has two important properties
\bea
 & &  ad \; cc'(t)  = ad\; c \circ ad\; c'(t), \nn \\
 & & ad \; c (t \otimes s) = {\displaystyle \sum_i\; (ad \; c_i (t) ) \otimes (ad \; c'_i (s))}. \nn 
\eea
From these properties, we can show that the adjoint action gives a 
representation of $ {\cal H} $
\beq
   ad \; [c, \; c'] (t) = [ad\; c,\; ad\; c'] (t).  \label{adrep}
\eeq
Tensor operators are defined as operators which form representation bases 
of a Hopf algebra under the adjoint action.
\begin{definition}
  Let $ T $ be a set of operators, and $ D(c)^{(j)} $ be a representation matrix 
of $ c \in {\cal H} $ with the highest weight $j$. The operators $ t_{jm} \in T $ 
are called rank $j$ tensor operators, if they satisfy the relations
\beq
  ad \; c(t_{jm}) = \sum_k\; D(c)^{(j)}_{km} t_{jk}.
\eeq
\end{definition}

  The adjoint action of $ X, Y $ and $ H $ is given by
\bea
 & & ad\; X(t_{jm}) = [X, \; t_{jm}], \nn \\
 & & ad\; Y(t_{jm}) = e^{-hX} [e^{hX}Y, \; t_{jm}] e^{-hX}, \label{adjoint} \\
 & & ad\; H(t_{jm}) = e^{-hX} [e^{hX}H, \; t_{jm}] e^{-hX}. \nn
\eea

\subsection{Some Examples of $ \slt $ Tensor Operators}

  In this section, we shall give explicit expressions of three kinds of 
$ \slt $ tensor operators. To show the existence of $ \slt $ tensor 
operators, it is enough to construct rank 1/2 tensor operators, since 
higher rank tensor operators can be obtained by decomposing a tensor 
product of some rank 1/2 tensor operators. This is due to  the 
fact that tensor operators are representation bases of $ \slt $ and 
we have a explicit formula for the $ \slt $ CGC.  

  The tensor operators 
given here are (1) rank 1/2 tensor operators in fermion realization of 
$ \slt $, (2) rank 1/2 tensor operators in boson realization of 
$ \slt $, (3) rank 1 tensor operators constructed by the generators of 
$ \slt $ themselves. The basic idea for (1) and (2) is quite simple.  
We realize $ \slt $ with the generators of $ sl(2) $
\beq 
 H = J_0, \qquad 
 X = {2 \over h} {\rm arctanh} \left({h \over 2} J_+\right), 
 \qquad 
 Y = \sqrt{1-\left({h \over 2}J_+\right)^2} J_- \sqrt{1-\left({h \over 2}J_+\right)^2},
   \label{real}
\eeq
where $ J_{\pm} $ and $ J_0 $ are generators of $ sl(2) $. This is obtained by 
solving (\ref{sltwo}) with respect to $ X $ and $ Y $ and regarding 
$ \{ Z_{\pm}, H \} $ as the generators of $ sl(2). $ 
Then we realize $ sl(2) $ in terms of fermions or bosons. 
We need representation matrices of $ X, Y $ and $ H $ for $ j=1/2 $ and $ 1 $ 
to find the rank 1/2 or 1 tensor operators.  
The representation matrices for $ j=1/2 $ read
\[
  X = \left(
  \begin{array}{cc}
  0 & 1 \\ 0 & 0
  \end{array}
  \right)_, \qquad
  Y = \left(
  \begin{array}{cc}
  0 & 0 \\ 1 & 0 
  \end{array}
  \right)_, \qquad
  H = \left(
  \begin{array}{cc}
  1 & 0 \\ 0 & -1
  \end{array}
  \right)_,
\] 
and for $ j=1 $
\bea
 & & 
 X = \left(
 \begin{array}{ccc}
 0 & \sqrt{2} & 0 \\ 0 & 0 & \sqrt{2} \\
 0 & 0 & 0
 \end{array}
 \right)_, \quad
 Y = \left(
 \begin{array}{ccc}
 0 & -h^2/2\sqrt{2} & 0 \\
 \sqrt{2} & 0 & -h^2/2\sqrt{2} \\
 0 & \sqrt{2} & 0
 \end{array}
 \right)_,
 \nn \\
 & &  H = diag(2, 0, -2). \nn
\eea
Note that the representation matrices for $ j=1/2 $ are same as the ones 
for $ h \rightarrow 1 $, however rank 1/2 tensor operators are 
nontrivial since the adjoint action has different form (see (\ref{adjoint})). 

  Let us first consider the fermion realization. We introduce two kinds of 
mutually anticommuting fermions
\[
 \{ a_i, \; a^{\dagger}_j \} = \delta_{ij}, \qquad 
 \{ a_i, \; a_j \} = \{ a^{\dagger}_i, \; a^{\dagger}_j \} = 0, \qquad
 i, j = 1, 2.
\]
These fermions realize $ sl(2) $ (the so-called fermion quasi spin 
formalism),
\beq
  J_+ = a^{\dagger}_1 a^{\dagger}_2, \qquad 
  J_- = a_2 a_1, \qquad
  J_0 = N_1 + N_2 - 1,  \label{quasispin}
\eeq
where $  N_i \equiv a^{\dagger}_i a_i $ is the number operator for $i$th 
fermion. This realization gives two representations of $ sl(2) $ and $ \slt $. 
One of them is the two dimensional irrep  whose representation space $ W^{(1/2)} $ 
has bases 
$ \ket{{1 \over 2} \; {1 \over 2}} = a_1^{\dagger} a_2^{\dagger} \ket{0} $ 
and $ \ket{{1 \over 2} \; -{1 \over 2}} = \ket{0} $, where $ \ket{0} $ 
denotes the fermion vacuum. The other is the trivial representation 
whose representation space $ W^{(0)} $ is spanned by $ a_1^{\dagger} \ket{0} $ 
or $ a_2^{\dagger} \ket{0}. $  The advantage of the fermions is that the adjoint 
action has a simpler form, since the nilpotency of fermions results 
$ X^2 = 0. $ We find two kinds of rank 1/2 tensor operators in this realization
\beq
  t_{1/2 \; 1/2} = -a_1^{\dagger}, \qquad
  t_{1/2 \; -1/2} = -a_2 + h (N_2-1) a_1^{\dagger}, \label{tenf1}
\eeq
and
\beq
   t_{1/2 \; 1/2} = a_2^{\dagger}, \qquad
    t_{1/2 \; -1/2} = -a_1 - h (N_1-1) a_2^{\dagger}. \label{tenf2}
\eeq
The straightforward computation shows that these satisfy the definition of 
rank 1/2 tensor operators. It is also easy to see that the action of both 
(\ref{tenf1}) and (\ref{tenf2}) on $ W^{(1/2)} $ results $ W^{(0)} $ and vice versa.

  Next we consider the boson realization. With two kinds of mutually commuting 
bosons
\[
  [b_i, \; b^{\dagger}_j] = \delta_{ij}, \qquad
  [b_i, \; b_j] = [b^{\dagger}_i, \; b^{\dagger}_j] = 0, \qquad 
  i, j = 1, 2,
\]
the Lie algebra $ sl(2) $ is realized as (the Jordan-Schwinger realization)
\beq
  J_+ = b^{\dagger}_1 b_2, \qquad
  J_- = b^{\dagger}_2 b_1, \qquad
  J_0 = N_1 - N_2,           \label{jordan}
\eeq
where $ N_i = b_i^{\dagger} b_i $ is the number operator for $i$th boson. 
We obtain any irrep of $ sl(2) $ and $ \slt $ in this realization. Let us 
denote the representation space for highest weight $j$ by $ W^{(j)} $, 
then the bases of $ W^{(j)} $ are given by 
\beq
    \ket{j m} = { (b_1^{\dagger})^{j+m} (b_2^{\dagger})^{j-m} \over 
                  \sqrt{(j+m)! (j-m)!} } \ket{0}, 
    \qquad
    m = -j, -j+1, \cdots, j,
\eeq
where $ \ket{0} $ denotes the boson vacuum. It is shown, by straightforward 
computation, that there exist two kinds of rank 1/2 tensor operators in 
this realization
\beq
  t_{1/2 \; 1/2} = \left( 1-{h \over 2}J_+\right)^{-1} b^{\dagger}_1,  \qquad
  t_{1/2 \; -1/2} =  \left( 1-{h \over 2}J_+\right) b^{\dagger}_2 + 
                       {h \over 2} (t_{1/2\; 1/2} - b^{\dagger}_1 J_0), 
  \label{tenb1}
\eeq
and
\beq
  t_{1/2 \; 1/2} = -\left( 1- {h \over 2}J_+\right)^{-1} b_2, \qquad
  t_{1/2 \; -1/2} =  \left( 1-{h \over 2}J_+\right) b_1 + 
                       {h \over 2} (t_{1/2\; 1/2} + b_2 J_0).
  \label{tenb2}
\eeq
The action of (\ref{tenb1}) on $ W^{(j)} $ reads
\bea
 & & t_{1/2\; 1/2} \ket{j m} = {\displaystyle 
       \sum_{n=0}^{j-m}\left({h \over 2}\right)^n \Gamma_n^{jm}
       \ket{j+{1\over 2} \;\;  m+{1 \over 2}+n}
       }, \label{action1} \\
 & & t_{1/2\; -1/2} \ket{j m} \nn \\
 & & = \sqrt{j-m+1} \ket{j+{1 \over 2}\; \; m-{1 \over 2}} 
       - {h \over 2} (j+m) \sqrt{j+m+1} \ket{j+{1 \over 2} \;\; m+{1 \over 2}} \nn \\
 & & +{\displaystyle 
       {h \over 2} \sum_{n=1}^{j-m}\; \left({h \over 2}\right)^n 
       \Gamma_n^{jm}
       \ket{j+{1 \over 2}\;\; m+{1 \over 2}+n}
       },\label{action2},
\eea
where
\[
  \Gamma_n^{jm} = \left\{
       { (j-m)! (j+m+n+1)! \over (j+m)! (j-m-n)!} \right\}^{1/2}_. 
\]
On the other hand, the action of (\ref{tenb2}) on $ W^{(j)} $ reads 
\bea
  & & t_{1/2\; 1/2} \ket{j m} = -{\displaystyle 
        \sum_{n=0}^{j-m-1}\left({h \over 2}\right)^n \Lambda_n^{jm}
        \ket{j-{1\over 2} \; \; m+{1 \over 2}+n}
        }, \label{action3} \\
  & &  t_{1/2\; -1/2} \ket{j m} \nn \\
  & & = \ket{j-{1 \over 2}\; \; m-{1 \over 2}} 
        - {h \over 2} \sqrt{j-m} (j-m-1) 
          \ket{j-{1 \over 2}\; \; m+{1 \over 2}+n} \nn \\
  & & - {\displaystyle {h \over 2} \sum_{n=1}^{j-m-1}\; 
        \left({h \over 2}\right)^n \Lambda_n^{jm} 
        \ket{j-{1\over 2}\; \; m+{1\over 2}+n}},
        \label{action4}
\eea
where
\[
   \Lambda_n^{jm} = {\displaystyle \left\{
        { (j-m)! (j+m+n)! \over (j+m)! (j-m-n-1)!}
        \right\}_.^{1/2} }
\]
Therefore, we see that the action of tensor operators (\ref{tenb1}) gives rise to a mapping 
$ W^{(j)} \rightarrow W^{(j+1/2)} $, while the tensor operators (\ref{tenb2}) 
$ W^{(j)} \rightarrow W^{(j-1/2)} $.

  The third example of tensor operators is constructed with 
the generators of $ \slt $ themselves. It is also straightforward to verify 
that the rank 1 tensor operators are given by
\bea
  & & 
  {\displaystyle t_{1\; 1} = -e^{hX} {\sinh hX \over h}}, \nn \\
  & & 
  {\displaystyle t_{1\; 0} = { e^{hX} H \over \sqrt{2}}}, \label{rank1} \\
  & & 
  {\displaystyle t_{1\; -1} = e^{-hX/2} Y e^{-hX/2} + {h \over 2} 
  e^{hX/2} H e^{hX/2} - {h \over 2} H^2}. \nn
\eea
These are combination of the $ \slt $ generators so that they can act 
on any irrep space and do not change the value of highest weight : 
$ t_{1\; m} : W^{(j)} \rightarrow W^{(j)} .$

 All the result given here is a natural analogue of  
$ sl(2) $, since they have the counterparts, which are well-known 
properties of $ sl(2) $ tensor operators,  in the 
limit of $ h \rightarrow 0. $ Therefore we have seen new similarities between 
the representation theories of $ \slt $ and $ sl(2) $.

%%%%%%%%%%%%%%%%%%%%%%%%%%%%%%%%%%%%%%%%%%%%%%%%%%%%%%
%
%    Wigner-Eckart's Theorem 
%
%%%%%%%%%%%%%%%%%%%%%%%%%%%%%%%%%%%%%%%%%%%%%%%%%%%%%%
%
\setcounter{equation}{0}
\section{Wigner-Eckart's Theorem}

  The results in the previous section enable us to consider an 
extension of the Wigner-Eckart's theorem to the Jordanian 
quantum algebra $ \slt. $ The purpose of this section is to 
show that the Wigner-Eckart's theorem can be extended to 
$ \slt. $ 

\begin{thm}
  Let $ T^{(j_1)} $ be a set of rank $ j_1 $ tensor operators, 
let $ W^{(j)} $  
be irrep space of  $ \;\slt $ with highest weight $ j $  
and suppose that 
$
    t_{j_1m_1}  \in T^{(j_1)} : W^{(j_2)} \rightarrow W^{(j)}. 
$
Then
\beq
 \bra{jm} t_{j_1 m_1} \ket{j_2 m_2} = I(j_1 j_2 j) \sum_{n_i = -j_i}^{j_i} \; 
 \alpha_{-m1, -m_2}^{-n_1, -n_2} C_{n_1,\; n_2,\; m}^{j_1,\; j_2,\; j},
 \label{wet}
\eeq
where $ I(j_1 j_2 j) $ is a constant independent of $ m_1, m_2 $ and $ m $. 
\end{thm}

\noindent
{\em Proof} : According to \cite{bt}, we consider an element 
$ t_{j_1 m_1} \otimes \ket{j_2  m_2} $ of $ T^{(j_1)} \otimes W^{(j_2)} $. 
Both $ T^{(j_1)} $ and $ W^{(j_2)} $ are representation space of $ \slt $ 
so that $ \Delta(c), \ c \in \slt $ acts on $ T^{(j_1)} \otimes W^{(j_2)} $.
For example,
\beq
  \Delta(H) t_{j_1 m_1} \otimes \ket{j_2  m_2} = 
  ad \; H(t_{j_1 m_1}) \otimes e^{hX} \ket{j_2  m_2} + ad\; e^{-hX} (t_{j_1 m_1}) 
  \otimes H \ket{j_2 m_2}.
  \label{induce}
\eeq
The LHS of $ \otimes $ is a tensor operator, since the adjoint action is 
a linear transformation on $ T^{(j_1)} $. Thus we can consider an action 
of the LHS of $ \otimes $ on the RHS : 
$ 
 t \otimes \ket{jm} \ \rightarrow \ t \ket{jm}.
$
This operation is called a ``contraction" in \cite{bt}. 
Noting that $ ad\; e^{-hX}(t_{j_1 m_1}) = e^{-hX} t_{j_1 m_1} e^{hX} $ 
and contractiong (\ref{induce}), we obtain $ H t_{j_1 m_1} \ket{j_2 m_2} $. 
Similar calculation shows that
\beq
  \Delta(c) t_{j_1 m_1} \otimes \ket{j_2  m_2} \ \rightarrow \ 
  c\; t_{j_1 m_1} \ket{j_2  m_2},
  \label{contract}
\eeq
where the arrow means that the RHS is a result of the contraction. 

  An intermediate vector on $ T^{(j_1)} \otimes W^{(j_2)} $ is given by
\beq
 \ket{(j_1 m_1) (j_2 m_2)} = \sum_{k_i = -j_i}^{j_i} \; \alpha_{k_1, k_2}^{m_1, m_2}\;
 t_{j_1 k_1} \otimes \ket{j_2 k_2}.
 \label{product}
\eeq
Because of (\ref{action}), the action of $ \Delta(Z_{\pm}) $ on the vector 
yields 
\bea
 \Delta(Z_{\pm}) \ket{(j_1m_1) (j_2m_2)}
 &=& \sqrt{(j_1\mp m_1) (j_1 \pm m_1 + 1)} \ket{(j_1\; m_1\pm 1) (j_2m_2)} 
 \nn \\
 &+& \sqrt{(j_2 \mp m_2) (j_2 \pm m_2 +1)} \ket{(j_1m_1) (j_2\; m_2\pm 1)}. 
 \label{actzpm}
\eea
Using (\ref{contract}), we obtain
\bea
  Z_{\pm} \ket{\varphi; m_1 m_2} &=& 
  \sqrt{(j_1 \mp m_1) (j_1 \pm m_1 + 1)} \ket{\varphi; m_1\pm 1 m_2}  \nn \\
  &+& 
   \sqrt{(j_2 \mp m_2) (j_2 \pm m_2 + 1)} \ket{\varphi; m_1 m_2\pm 1},
 \label{rec1}
\eea
where $ \ket{\varphi; m_1 m_2} $ is the vector obtained from (\ref{product}) 
by a contraction.
\[
 \ket{\varphi; m_1 m_2} =
 \sum_{k_i = -j_i}^{j_i} \; \alpha_{k_1, k_2}^{m_1, m_2}\;
 t_{j_1 k_1} \ket{j_2 k_2}.
\]
The inner product of $ \ket{j \; m\pm 1} $ and (\ref{rec1}) gives the 
recurrence relations for $ \braket{jm}{\varphi; m_1 m_2} $ 
\bea
 & & \sqrt{(j\mp m) (j \pm m + 1)} \braket{jm}{\varphi; m_1 m_2} \nn \\
 & & \hspace{1cm}
       = \sqrt{(j_1\mp m_1) (j_1 \pm m_1 + 1)} 
       \braket{j \; m\pm 1}{\varphi; m_1\pm 1 \; m_2} \label{rec2}  \\
 & & \hspace{1cm}
       + \sqrt{(j_2 \mp m_2) (j_2 \pm m_2+1)} 
       \braket{j\; m\pm 1}{\varphi; m_1\; m_2 \pm 1}. \nn
\eea
The recurrence relations (\ref{rec2}) are the same as the ones for $ sl(2) $ CGC, 
therefore the quantity $ \braket{jm}{\varphi; m_1 m_2} $ must be 
proportional to $ sl(2) $ CGC. Denoting the proportional coefficient by 
$ I(j_1 j_2 j) $
\beq
  \braket{jm}{\varphi; m_1 m_2} = \sum_{k_i = -j_i}^{j_i}\; 
  \alpha_{k_1, k_2}^{m_1, m_2}\; 
  \bra{jm} t_{j_1 k_1} \ket{j_2 k_2} 
  =
  C_{m_1,\; m_2,\; m}^{j_1,\; j_2,\; j} I(j_1j_2j).
  \label{CG}
\eeq
This relation can be solved with respect to $ \bra{jm} t_{j_1 k_1} \ket{j_2 k_2}  $, and  
the Wigner-Eckart's theorem (\ref{wet}) has been proved.

\noindent
{\em Remark} : From (\ref{bra}) and the fact that the $ sl(2) $ CGC for bra 
and ket vectors are equal, we see that the quantity appearing in the right hand 
side of (\ref{wet}) is the $ \slt $ CGC for bra vectors. This is a general 
property of the Wigner-Eckart's theorem \cite{naru2}.

\end{document}